\documentclass[12pt,a4paper]{amsart}

\usepackage{fancyhdr}
\usepackage{amsmath, amsfonts, amssymb, amsthm}
\usepackage[english]{babel}
\usepackage[all]{xy}
\usepackage{graphicx}
\usepackage[T1]{fontenc}
\usepackage[utf8]{inputenc}
\usepackage{color}
\usepackage[a4paper,top=3cm,bottom=3cm,left=3cm,right=3cm]{geometry}
\usepackage{tikz}
\usetikzlibrary{matrix, arrows, decorations.pathmorphing}
\usepackage{plain}
\usepackage{url}

\title{\textbf{$N$-fiber-full modules}}
\author{Hongmiao Yu}

\newtheoremstyle{break}%
{}{}%
{\slshape}{}%
{\bfseries}{.}%
{5pt}{}

\theoremstyle{definition}
\newtheorem{definition}{Definition}[section]

\newtheorem{example}[definition]{Example}

\newtheorem{notation}[definition]{Notation}
\newtheorem{notations}[definition]{Notations}
\newtheorem{remark}[definition]{Remark}

\theoremstyle{break}
\newtheorem{lemma}[definition]{Lemma}
\newtheorem{proposition}[definition]{Proposition}

\newtheorem{theorem}[definition]{Theorem}
\newtheorem{corollary}[definition]{Corollary}

\newtheorem*{maintheorem*}{Main\ Theorem(cf. Theorem\til\ref{main})}
\newtheorem*{theorem*}{Theorem}
\newtheorem*{theorem3.5*}{Theorem(cf. Theorem\til\ref{thm35})}

\def\til{~}

\newcommand\tbs[1]{\textsl{\textbf{#1}}}  

\def\N{\mathbb N}

\def\Z{\mathbb Z}

\def\mm{\mathfrak m}  

\def\qq{\mathfrak q}

\def\H{H} 

\def\la{\longrightarrow}
\def\map{\longmapsto}
\def\:{\colon}

\def\Hom{\mathrm{Hom}}

\def\Ker{\mathrm{Ker}}
\def\Im{\mathrm{Im}}

\def\altezza{\mathrm{ht}}   

\def\Ext{\mathrm{Ext}}

\def\supp{\mathrm{Supp}}

\def\height{\mathrm{ht}}

\def\ass{\mathrm{Ass}}

\def\cano{\mathrm\omega}

\def\min{\mathrm{Min}}

\def\se{\subseteq}

\def\cl{\overline}

\def\iso{\cong}

\def\init{\mathrm{in}}

\def\initi{\mathrm{init}}

\def\proj{\mathrm{Proj}}

\def\satu{\mathrm{sat}}

\def\cocoa{\mbox{\rm 
   C\kern-.13em o\kern-.07 em C\kern-.13em o\kern-.15em A}}

\newcommand\fg[0]{finitely generated{}}

\newcommand\Nff[0]{$N$-fiber-full{}}

\newcommand\sat[1]{{#1}^\satu}

\begin{document}
\begin{abstract}
Let $A$ be  a Noetherian flat $K[t]$-algebra, $h$  an integer and let $N$ be a graded $K[t]$-module, we introduce and study ``$N$-fiber-full up to $h$'' $A$-modules. We prove that an $A$-module $M$ is \Nff\ up to $h$ if and only if $\Ext^i_A(M, N)$ is flat over $K[t]$ for all $i\le h-1$. And we show some applications of this result extending the recent result on square-free Gr\"obner degenerations by Conca  and Varbaro.
\end{abstract}

\maketitle

\section{Introduction}
Throughout this paper, $A$ is a Noetherian flat $K[t]$-algebra, $M$ and $N$ are \fg\ $A$-modules which are flat over $K[t]$, and all of $A$, $M$ and $N$ are graded $K[t]$-modules. \\
At the ``CIME-CIRM Course on Recent Developments in Commutative Algebra'' conference in 2019, Matteo Varbaro introduced the notion of ``fiber-full modules''  providing a new proof of the main result of \cite{CV}. The starting point of this paper is to find  some possible generalizations  of this concept. We recall that $M$ is a \tbs{fiber-full} $A$-module if, for any $m\in\N_{>0}$, the natural projection $M/t^mM\la M/tM$ induces injective maps $\Ext^i_A(M/tM, A)\la\Ext^i_A(M/t^mM, A)$  for all $i\in \Z$. And one of the ``most important''  properties of fiber-full modules is that 
being fiber-full is related to the flatness over $K[t]$ of every $\Ext^i_A(M, A)$ \cite{levico}.
We notice that if $N$ is a graded $K[t]$-module, the induced map $\Ext^0_A(M/tM, N)\la\Ext^0_A(M/t^mM, N)$ is always injective. Motivated by this, let us introduce the ``$N$-fiber-full up to $h$'' modules: 
\begin{definition}
Let $h$ be an integer. We say that $M$ is \tbs{\Nff\ up to $h$} as an $A$-module if, for any $m\in\N_{>0}$, the natural projection $M/t^mM\la M/tM$ induces injective maps $\Ext^i_A(M/tM, N)\la\Ext^i_A(M/t^mM, N)$  for all $i\le h$.
\end{definition}
An important question is: if $M$ is \Nff\ up to $h$ as an $A$-module, can we obtain the flatness of some $\Ext^i_A(M, N)$? The main theorem of this paper goes in this direction.
\begin{maintheorem*}
Let $h$ be an integer. $M$ is \Nff\ up to $h$ as an $A$-module if and only if $\Ext^i_A(M, N)$ is flat over $K[t]$ for all $i\le h-1$.
\end{maintheorem*}

To prove this theorem, the most difficult part  is the proof of Lemma\til\ref{flatness} below, concerning the equivalence between two properties 1) and 2). At the beginning, I thought one could prove it using a way similar to the proof of \textcolor{black}{Lemma $3.5$} in \cite{levico}, actually the proof of implication $1\Rightarrow 2)$ is done in this way. But the proof of the other implication of this lemma is the hardest point of this paper: the whole section $2$ is to prove Lemma\til\ref{flatness} and the main theorem. In section $3$, we will talk about some applications of $N$-fiber-full modules. A main consequence, as we will see, is that the notion ``$N$-fiber-full up to $h$'' allows us to infer interesting results whenever the special fiber $M/tM$ has ``nice'' properties after removing primary components of big height. For example, 
\begin{theorem3.5*}
Let $S$  be the polynomial ring $K[X_1,\dots, X_n]$ over a field $K$, let  $I\se S$ be a homogeneous ideal. Fixing a monomial order on $S$, we denote by $\init(I)$ the initial ideal of $I$ with respect to this monomial order. If $I$ is such that $\sat {\init(I)}$ is square-free, then $$\dim_K\H_\mm^i(S/I)_j=\dim_K\H_\mm^i(S/\init(I))_j$$ for all $i\ge 2$ and for all $j\in\Z$.
\end{theorem3.5*}
Equivalently,
\begin{theorem*}
Let $\mathbb P^n$ be the n-dimensional projective space over a field $K$, let $X\se \mathbb P^n$ be a projective scheme and let $I$ be a homogeneous ideal of $S=K[X_0, \dots, X_n]$ such that $X=\proj(S/I)$. Fix a monomial order $<$ on $S$ and 
assume that 
$\init(X)=\proj(S/\init(I))$ is reduced, where $\init(I)$ is the initial ideal of $I$ with respect to $<$. Then
$$\dim_K\H^i(X, \mathcal O_X(j))=\dim_K\H^i(\init(X), \mathcal O_{\init(X)}(j))$$
for all $i>0$ and for all $j\in\Z$.
\end{theorem*}

This theorem has already been announced by Varbaro in his paper \cite{umi} (Theorem $4.4$) without writing the complete proof.
By a private communication, Varbaro told me that he realized later that it was not clear to him how to extend the proof given in \cite{levico}
for the same reasons explained above. But we will show that we can complete the proof of a very general version of this theorem using    Theorem\til\ref{main} of this paper.

\section{definition and properties of $N$-fiber-full modules}


In this section we prove some basic properties  of $N$-fiber-full modules, in order to show the equivalence between being $N$-fiber-full up to $h$ and the flatness of  modules $\Ext^i_A(M, N)$ with $i\le h-1$.

\begin{remark}
If $A$ is a Cohen-Macaulay complete local ring, $N$ is the canonical module of $A$ and $M=A/I$ with $I$ an ideal of $A$, then $M$ is \Nff\ up to $\dim A$ is equivalent to saying that $t$ is a surjective element in $M$ in the sense of \cite{DDSM} section $3.1$.\\
If $A$ is a local ring and $N=A$, then $M$ is \Nff\ up to $h$ for each $h\in\N$ is equivalent to saying that $M$ is fiber-full defined as in \cite{levico} \textcolor{black}{Definition $3.8$}.
\end{remark}

One implication of the main theorem is not difficult to prove,  but it is very useful because we need it to prove the other one.

\begin{theorem}\label{leftarrow}
Let $h$ be an integer and let $\Ext^i_A(M, N)$ be flat over $K[t]$ for all $i\le h-1$. Then $M$ is \Nff\ up to $h$ as an $A$-module.
\proof
If $\Ext^i_A(M, N)$ is flat over $K[t]$ for all $i\le h-1$, then $$\Ext^i_A(M, N)\xrightarrow{\ \cdot t^{m-1}\ }\Ext^i_A(M, N)$$ is injective for all $i\le h-1$ and for all $m\in\N_{>0}$, hence $$\Ext^{i-1}_A(M, N)\xrightarrow{\ \cdot t^{m-1}\ }\Ext^{i-1}_A(M, N)$$ is injective for all $i\le h$ and for all $m\in\N_{>0}$. The commutative diagram of $A$-modules with exact rows
$$\xymatrix@=3em{
0\ar[r] &M \ar[r]^{\cdot t} & M\ar[r] &M/tM\ar[r]  &0\\
0\ar[r] &M \ar[r]^{\cdot t^m}\ar[u]^{t^{m-1}}& M\ar[r]\ar@{=}[u] &M/t^mM\ar[r]\ar[u]&0
}$$
yields the following commutative diagram of $A$-modules with exact rows
$$\xymatrix@=3em{
\Ext^{i-1}_A(M,N) \ar[r]^{\cdot t}\ar@{=}[d]& \Ext^{i-1}_A(M,N)\ar[r] \ar[d]^{t^{m-1}} &\Ext^i_A(M/tM,N)\ar[r]\ar[d]  &\Ext^i_A(M,N)\ar@{=}[d]\\
\Ext^{i-1}_A(M,N) \ar[r]^{\cdot t^m}& \Ext^{i-1}_A(M,N)\ar[r] &\Ext^i_A(M/t^mM,N)\ar[r]&\Ext^i_A(M,N).
}$$
Since $$\Ext^{i-1}_A(M, N)\xrightarrow{\ \cdot t^{m-1}\ }\Ext^{i-1}_A(M, N)$$ is injective for all $i\le h$, by the proof of the Five Lemma we have that 
$$\Ext^i_A(M/tM, N)\la\Ext^i_A(M/t^mM, N)$$ is injective  for all $i\le h$.
\endproof
\end{theorem}
 
\begin{remark}
The map $M/t^kM\xrightarrow{\cdot t^{l-k}} M/t^lM$ is injective for each $k,l\in\N_{>0}$ such that $l\ge k$.
\proof
If $x\in M$ is such that $\cl{t^{l-k}x}=0$ in $M/t^lM$, then there exists $y\in M$ such that $t^ly=t^{l-k}x$, hence $t^{l-k}(x-t^ky)=0$. Since $M$ is flat over $K[t]$, we have $x=t^ky$, it follows that $\cl x=0$ in $M/t^kM$.\endproof
\end{remark}

\begin{lemma}\label{MAQUY}
The following are equivalent:
\begin{enumerate}
\item[$1$)] $M$ is \Nff\ up to $h$.
\item[$2$)] For each $k,l\in\N_{>0}$ such that $l\ge k$, the short exact sequence 
$$0\la M/t^kM\xrightarrow{\cdot t^{l-k}} M/t^lM\la M/t^{l-k}M\la0$$
induces a short exact sequence
$$0\la\Ext^i_A(M/t^{l-k}M,N)\la\Ext^i_A(M/t^lM, N)\xrightarrow{\ f^i_{k,l}\ }\Ext^i_A(M/t^kM, N)\la 0 $$ for all $i\le h-1$.
\item[$3$)] For each $k, l\in\N_{>0}$ such that $l\ge k$, the natural projection $M/t^lM\la M/t^kM$ induces injective maps $\Ext^i_A(M/t^kM, N)\la\Ext^i_A(M/t^lM, N)$  for all $i\le h$.
\end{enumerate}
Furthermore, up to the identifications $$t^{l-k}\Ext^i_A(M/t^lM, N)\iso\Ext^i_A(M/t^kM, N)\text{\quad for each $i\le h-1$,}$$ the map $f^i_{k,l}$ in $2)$ corresponds to the surjective map
$$\Ext^i_A(M/t^lM, N)\xrightarrow{\cdot t^{l-k}}t^{l-k}\Ext^i_A(M/t^lM, N)$$
for each $i\le h-1$.
\proof
In the case $k=l$, $2$) and $3$) are always true, so we suppose $k<l$.
\begin{enumerate}
\item[$1\Rightarrow 2$)] We consider the short exact sequence 
$$0\la M/t^kM\stackrel{\cdot t}\la M/t^{k+1}M\la M/tM\la0,$$
it induces a long exact sequence
$$\dots\xrightarrow{\ f^{i-1}_{k,k+1}\ }\!\Ext^{i-1}_A(M/t^kM, N) \!\la\!\Ext^{i}_A(M/tM,N)\!\la\!\Ext^{i}_A(M/t^{k+1}M, N)\!\la\!\dots\!. $$
Since $M$ is \Nff\ up to $h$, $\Ext^{i}_A(M/tM,N)\la\Ext^{i}_A(M/t^{k+1}M, N)$ is injective for each $i\le h$, it follows that $f^{j}_{k,k+1}$ is surjective for each $j\le h-1$. Hence $f^{j}_{k,l}= f^{j}_{k,k+1}\circ\dots\circ f^{j}_{l-1,l}$ is surjective for each $j\le h-1$. 
Therefore $$0\la\Ext^j_A(M/t^{l-k}M,N)\la\Ext^j_A(M/t^lM, N)\xrightarrow{\ f^j_{k,l}\ }\Ext^j_A(M/t^kM, N)\la 0 $$
is exact for all $j\le h-1$.
\item[$2\Rightarrow 3$)] The short exact sequence 
$$0\la M/t^{l-k}M\xrightarrow{\cdot t^{k}} M/t^lM\la M/t^{k}M\la0$$
induces a long exact sequence
$$\!\dots\!\xrightarrow{\ f^{i}_{l-k,l}\ }\!\Ext^{i}_A(M/t^{l-k}M, N) \!\la\!\Ext^{i+1}_A(M/t^kM,N)\!\la\!\Ext^{i+1}_A(M/t^{l}M,N)\!\la\!\dots\!. $$
By $2$) $f^{i}_{l-k,l}$ is surjective for all $i\le h-1$, it follows that $$\Ext^{j}_A(M/t^kM,N)\la\Ext^{j}_A(M/t^{l}M, N)$$ is injective for all $j\le h$.
\item[$3\Rightarrow 1$)] Take $k=1$.
\end{enumerate}
Furthermore, we observe that the endomorphism $$M/t^lM\xrightarrow{\cdot t^{l-k}}  M/t^lM$$ corresponds to the composition of maps
$$M/t^lM\stackrel{p}\la M/t^kM\xrightarrow{\cdot t^{l-k}}  M/t^lM$$ where $p$ is the natural projection, therefore the following diagram
$$\xymatrix@=3em{
\Ext^i_A(M/t^lM, N)\ar[rr]^{\cdot t^{l-k}}\ar[rd]_{f^i_{k,l}}&&\Ext^i_A(M/t^lM, N)\\
&\Ext^i_A(M/t^kM, N)\ar[ru]&
}$$
is commutative.  By $3$) $\Ext^i_A(M/t^kM, N)\la\Ext^i_A(M/t^lM, N)$ is injective for all $i\le h-1$, and by $2$) $f^i_{k,l}$ is surjective for all $i\le h-1$. We have that for each $i\le h-1$, 
$$\Ext^i_A(M/t^kM, N)\iso t^{l-k}\Ext^i_A(M/t^lM, N)$$ is a submodule of $\Ext^i_A(M/t^lM, N)$, and  $f^i_{k,l}$ corresponds to the surjective map
$$\Ext^i_A(M/t^lM, N)\xrightarrow{\cdot t^{l-k}}t^{l-k}\Ext^i_A(M/t^lM, N).$$
\endproof
\end{lemma}

The next lemma is crucial to show Theorem\til\ref{main}. The most difficult implication $2\Rightarrow 1)$ surprisingly uses Theorem\til\ref{leftarrow}.

\begin{lemma}\label{flatness}
Let $h$ be an integer. The following are equivalent:
\begin{enumerate}
\item[$1$)] $\Ext_{A}^i(M, N)$ is flat over $K[t]$ for all $i\le h$.
\item[$2$)] $\Ext_{A/t^mA}^i(M/t^mM, N/t^mN)$ is flat over $K[t]/(t^m)$ for all $m\in\N_{>0}$ and for all $i\le h-1$.
\end{enumerate}
\proof  
$1\Rightarrow 2)$: Since $N$ is flat over $K[t]$, there is a short exact sequence 
$$0\la N\stackrel{\cdot t^m}\la N\la N/t^mN\la 0.$$
Consider the induced long exact sequence of $\Ext^i_A(M,-):$
$$\xymatrix@=1.5em{
\dots\ar[r]&\Ext^i_A(M,N) \ar[r]^{\cdot t^m}&\Ext^i_A(M,N)\ar[r]&\Ext^i_A(M,N/t^mN)\ar[r]&\\
&\Ext^{i+1}_A(M,N) \ar[r]^{\cdot t^m}&\Ext^{i+1}_A(M,N)\ar[r]&\dots.\quad\quad\quad\quad
}$$
Since $\Ext_{A}^k(M, N)$ is flat over $K[t]$ for all $k\le h$,  $t^m$ is an $\Ext_{A}^k(M, N)$-regular element for all $m\in\Z_+$, and so we obtain a short exact sequence
$$0\la\Ext^i_A(M,N)\stackrel{\cdot t^m}\la\Ext^i_A(M,N)\la\Ext^i_A(M,N/t^mN)\la0$$
for all $i\le h-1$. It follows that
$$\Ext^i_A(M,N/t^mN)\iso\frac{\Ext^i_A(M,N)}{t^m\Ext^i_A(M,N)}$$
for all $i\le h-1$.
Furthermore, using again $1$) we have that $\Ext^i_A(M,N)/t^m\Ext^i_A(M,N)$ is flat over $k[t]/(t^m)$ for $i\le h$ (see \cite{matsumura} Section $7$). Therefore,
$$\Ext_{A/t^mA}^i(M/t^mM, N/t^mN)\iso\Ext^i_A(M,N/t^mN)$$ is flat over $k[t]/(t^m)$ for all $i\le h-1$.
\\
$1\Leftarrow 2)$: We use induction on $h\ge0.$ If $h=0$, we consider the long exact sequence $$0\la\Hom_A(M/tM,N)\la \Hom_A(M,N)\stackrel{\cdot t}\la\Hom_A(M,N)\la\dots,$$ 
induced by the short exact sequence
$0\la M\stackrel{\cdot t}\la M\la M/tM\la 0.$
Since $N$ is flat over $K[t]$, we have $\Hom_A(M/tM,N)=0$, and it follows that the map $\Hom_A(M,N)\stackrel{\cdot t}\la\Hom_A(M,N)$ is injective. Hence  $1)$ is true. \\
If $h\ge 1$, we suppose that $\Ext_{A/t^mA}^i(M/t^mM, N/t^mN)$ is flat over $K[t]/(t^m)$ for all $m\in\N_{>0}$ and for all $i\le h-1$. 
By the inductive hypothesis $\Ext_{A}^i(M, N)$ is flat over $K[t]$ for all $i\le h-1$. Furthermore, $M$ is \Nff\ up to $h$ by  Theorem\til\ref{leftarrow}. We prove that $\Ext_{A}^h(M, N)$ is flat over $K[t]$.
By contradiction, suppose that  there exists $x\in\Ext^h_A(M,N)$, $x\not=0$ and there exists $k\in\N$ such that $t^kx=0$. Let $l$ be an integer such that $l>k$. The commutative diagram 
$$\xymatrix@=3em{
0\ar[r] &M \ar[r]^{\cdot t^l} & M\ar[r] &M/t^lM\ar[r]  &0\\
0\ar[r] &M \ar[r]^{\cdot t^{l-k}}\ar@{=}[u]& M\ar[r]\ar[u]_{\cdot t^{k}} &M/t^{l-k}M\ar[r]\ar[u]_{\cdot t^{k}}&0
}$$
yields the following commutative diagram of $A$-modules with exact rows
$$\xymatrix@=3em{
\Ext^{h-1}_A(M,N) \ar[r]^\psi\ar@{=}[d]& \Ext^{h}_A(M/t^lM,N)\ar[r]^\phi \ar[d]^{f^h_{l-k,l}} &\Ext^h_A(M,N)\ar[r]^{\cdot t^l}\ar[d]^{\cdot t^k} &\Ext^h_A(M,N)\ar@{=}[d]\\
\Ext^{h-1}_A(M,N) \ar[r]^\eta& \Ext^{h}_A(M/t^{l-k}M,N)\ar[r]^\zeta &\Ext^h_A(M,N)\ar[r]^{\cdot t^{l-k}}&\Ext^h_A(M,N).
}$$
Since $t^kx=0$, $l>k$, we have $t^lx=0$, hence there exists $y\in\Ext^{h}_A(M/t^lM,N)$, $y\not=0$ such that $\phi(y)=x$. We may suppose that $f^h_{l-k,l}(y)=0$. Indeed, if $f^h_{l-k,l}(y)\not=0$, since $\zeta(f^h_{l-k,l}(y))=t^k(\phi(y))=t^kx=0$, there exists $z\in \Ext^{h-1}_A(M,N)$ such that $\eta(z)=f^h_{l-k,l}(y)$. Set $y'=y-\psi(z)$. We have $$\phi(y')=\phi(y)-\phi(\psi(z))=x$$ and $$f^h_{l-k,l}(y')=f^h_{l-k,l}(y)-f^h_{l-k,l}(\psi(z))=f^h_{l-k,l}(y)-\eta(z)=0.$$ 
Hence, using Lemma\til\ref{MAQUY}, we may suppose that $t^ky=0$ in $\Ext^{h}_A(M/t^lM,N)$.
Furthermore, since  $\Ext^{h}_A(M/t^lM,N)\iso\Ext^{h-1}_{A/t^lA}(M/t^lM,N/t^lN)$ (see \cite{CMrings} Lemma $3.1.16$), there exists $[\varphi]\in\Ext^{h-1}_{A/t^lA}(M/t^lM,N/t^lN)$, $[\varphi]\not=0$ such that $t^k[\varphi]=0$.\\
Since $M$ is \Nff\ up to $h$, the natural map $$\Ext^h_A(M/t^lM,N)\la\Ext^h_A(M/t^mM,N)$$ is injective for each $m\ge l$ by Lemma\til\ref{MAQUY}. Hence $[\varphi]\in\Ext^{h-1}_{A/t^mA}(M/t^mM,N/t^mN)$ for each $m\ge l$ and $t^k[\varphi]=0$ in $\Ext^{h-1}_{A/t^mA}(M/t^mM,N/t^mN)$ for each $m\ge l$. Now it is enough to find a positive integer $m\ge l$ such that $[\varphi]\not\in t^{m-k}\Ext^{h-1}_{A/t^mA}(M/t^mM,N/t^mN)$.\\
Let us take an $A$-free resolution $F_\bullet$ of $M$, and let $(G^\bullet, \partial^\bullet)$ be the complex $$\Hom_A(F_\bullet, N),$$ so that $\Ext^i_A(M,N)$ is the $i$-th cohomology module of $G^\bullet$. Since $M$ and $A$ are flat over $K[t]$, $F_\bullet/t^mF_\bullet$ is an $A/t^mA$-free resolution of $M/t^mM$ for all $m\in\N_{>0}$. Let $(\cl {G^\bullet}, \cl{\partial^\bullet})$ denote the complex $$\Hom_{A/t^mA}(F_\bullet/t^mF_\bullet, N/t^mN),$$ so that $\Ext^i_{A/t^mA}(M/t^mM,N/t^mN)$ is the $i$-th cohomology module of $\cl{G^\bullet}$ and $\pi^\bullet$ the natural map of complexes from $G^\bullet$ to  $\cl{G^\bullet}$. Since $F_{h-1}$ is free and \fg, there exists $\beta_{h-1}\in\N$ such that $F_{h-1}=A^{\beta_{h-1}}$, hence
$$G^{h-1}=\Hom_A(F_{h-1}, N)=\Hom_A(A^{\beta_{h-1}}, N)\iso N^{\beta_{h-1}}$$
and 
$$\cl {G^{h-1}}=\Hom_{A/t^mA}(F_{h-1}/t^mF_{h-1}, N/t^mN)\iso (N/t^mN)^{\beta_{h-1}},$$
therefore $\pi^{h-1}: G^{h-1}\la \cl {G^{h-1}}$ is surjective. It follows that there exists $\delta\in G^{h-1}$ such that $\pi^{h-1}(\delta)=\varphi$.
$$\xymatrix@=4em{
G^{h-1}\ar[r]^{\partial^{h-1}}\ar[d]^{\pi^{h-1}}&G^h\ar[r]^{\partial^{h}}\ar[d]^{\pi^{h}}&G^{h+1}\\
\cl {G^{h-1}}\ar[r]^{\cl{\partial^{h-1}}}&\cl {G^h}&.
}$$
Set $w=\partial^{h-1}(\delta)\in G^h$. Since $\cl{\partial^{h-1}}(\varphi)=0$, then 
\begin{align*}
\pi^h(w)=0\in \cl {G^h}&=\Hom_{A/t^mA}(F_h/t^mF_h,N/t^mN)\\
&\iso\Hom_{A}(F_h,N)/t^m \Hom_{A}(F_h,N).
\end{align*}
So there exists $w'\in G^h$ such that $w=t^mw'$. Hence $w\in t^mG^h$ for all $m\gg0$. Using Krull's intersection theorem we have that $\partial^{h-1}(\delta)=w=0$, thus $[\delta]\in\Ext^{h-1}_{A}(M,N)$. \\
Since $[\varphi]\not=0$, $\varphi\not\in\Im(\cl {\partial^{h-2}})$.
If $\varphi=\pi^{h-1}(\delta)\in t^{m-k}\ker(\cl{\partial^{h-1}})$, 
$$\delta\in t^{m-k}\ker(\partial^{h-1})+t^mG^{h-1}\se t^{m-k}G^{h-1}.$$
But $\delta\not=0$, so $\varphi\not\in t^{m-k}\ker(\cl{\partial^{h-1}})$ for all $m\gg0$ by Krull's intersection theorem.\\
Therefore,  there exists $m$ such that $$t^k[\varphi]=0$$ and $$[\varphi]\not\in t^{m-k}\Ext^{h-1}_{A/t^mA}(M/t^mM,N/t^mN),$$ but this contradicts the fact that $\Ext_{A/t^mA}^{h-1}(M/t^mM, N/t^mN)$ is flat over $K[t]/(t^m)$.
\endproof
\end{lemma}

%
%
%

\begin{notations}
We introduce some notations which are useful in the  following
results.
\begin{itemize}
\item $A_m=A/t^mA$ for each $m\in\Z_+$.
\item $M_m=M/t^mM$  for each $m\in\Z_+$.
\item $N_m=N/t^mN$  for each $m\in\Z_+$.
\item $\iota_j: t^{j+1}M_m\la t^j M_m$ the natural inclusion  for each $m\in\Z_+$ and for each $j\in\N$.
\item $\mu_j: t^j M_m\la t^{m-1} M_m$ the multiplication by $t^{m-1-j}$ for each $m\in\Z_+$ and for each $j\in\N$, $j\le m-1$.
\item $E^i_m(-)$ the contravariant functor $\Ext^i_{A_m}(-, N_m)$ for each $i\in\N$.
\end{itemize}
\end{notations}

\begin{remark}\label{rm1}
Suppose that $k\in\Z_+$. Since $t^m$ is an $A$- and $N$-regular element for all $m\in\Z_+$ and $t^mM_k=0$ for all $m\ge k$, we have that
$$\Ext^{i+1}_A(M_k,N)\iso\Ext^{i}_{A/t^m}(M_k,N/t^mN)=E^i_m(M_k)$$
for all $i\in\N$ and for all $m\ge k$ by Lemma $3.1.16$ in \cite{CMrings}. In particular, 
$$\Ext^{i+1}_A(M_k,N)\iso E^i_k(M_k)$$ for all $i\in\N$.
Therefore, 
$$E^i_k(M_k)\iso \Ext^{i+1}_A(M_l,N)\iso E^i_m(M_k)$$
for all $i\in\N$ and for all $k\le m$.
\end{remark}

\begin{remark}\label{rm2}
The following map
\begin{eqnarray*}
M_j&\la& t^{m-j}M_m\\
\cl m&\map& [t^{m-j}m]
\end{eqnarray*}
is an isomorphism for all $j\le m$. 
\proof
If $j=m$, it is trivial. We suppose that $j<m$. It is clear that  the above map is surjective. If $t^{m-j}m\in t^mM$, then there exists $n\in M$ such that $t^{m-j}(m-t^{j}n)=t^{m-j}m-t^mn=0$. Since $t^{m-j}$ is an $M$-regular element for all $j<m$, $m=t^{j}n\in t^jM$. It follows that the above map is injective. \endproof
\end{remark}

\begin{remark}\label{rm3}
Let $h$ be an integer and let $j\le m-1$. If $M$ is \Nff\ up to $h$, then the following short exact sequence 
$$0\la t^{j+1}M_m\stackrel{\iota_j}\longrightarrow t^jM_m\stackrel{\mu_j}\longrightarrow t^{m-1}M_m\la 0$$
yields a short exact sequence
$$0\la E^i_m(t^{m-1}M_m)\xrightarrow{\ E^i_m(\mu_j)\ }E^i_m(t^jM_m)\xrightarrow{E^i_m(\iota_j)}E^i_m(t^{j+1}M_m)\la 0$$
for all $i\le h-2$.
\proof
It is easy to check that $0\la t^{j+1}M_m\stackrel{\iota_j}\longrightarrow t^jM_m\stackrel{\mu_j}\longrightarrow t^{m-1}M_m\la 0$ is a short exact sequence and, using the identifications in  Remark\til\ref{rm2}, the following diagram is commutative
$$\xymatrix@=3em{
0\ar[r] &t^{j+1}M_m\ar[r]^{\iota_j} &t^jM_m\ar[r]^{\mu_j}& t^{m-1}M_m\ar[r]&0\\
0\ar[r] &t^{j+1}M_m\ar[r]\ar@{=}[u]&M_{m-j}\ar[r]\ar[u]_\iso& M_1\ar[r]\ar[u]_\iso&0.
}$$
 Applying $E^i_m(-)$ on the short exact sequence and the above diagram we obtain a long exact sequence \begin{eqnarray*}
\dots\la E^i_m(t^{m-1}M_m)&\xrightarrow{E^i_m(\mu_j)}&E^i_m(t^jM_m)
\xrightarrow{E^i_m(\iota_j)}E^i_m(t^{j+1}M_m)\la E^{i+1}_m(t^{m-1}M_m)\\
&\xrightarrow{E^{i+1}_m(\mu_j)}&E^{i+1}_m(t^jM_m)\la \dots
\end{eqnarray*}
and a commutative diagram 
$$\xymatrix@=3em{
E^i_m(t^{m-1}M_m)\ar[d]^\iso\ar[r]^{E^i_m(\mu_j)}&E^i_m(t^jM_m)\ar[d]^\iso\\
E^i_m(M_1)\ar[r]&E^i_m(M_{m-j})
}$$
for all $i\in\N$.
Since $M$ is \Nff\ up to $h$, $\Ext^j_A(M/tM, N)\la\Ext^j_A(M/t^mM, N)$ is injective for all $j\le h$, and so  $E^i_m(M_1)\la E^i_m(M_{m-j})$ is injective for all $i\le h-1$ by Remark\til\ref{rm1}. Hence $E^i_m(\mu_j)$ is injective for all $i\le h-1$ by the above diagram. It follows that 
$$0\la E^i_m(t^{m-1}M_m)\xrightarrow{\ E^i_m(\mu_j)\ }E^i_m(t^jM_m)\xrightarrow{E^i_m(\iota_j)}E^i_m(t^{j+1}M_m)\la 0$$ is exact for all $i\le h-2$.
\endproof
\end{remark}

\begin{theorem}\label{main}
Let $h$ be an integer. 
$M$ is \Nff\ up to $h$ as an $A$-module if and only if $\Ext^i_A(M, N)$ is flat over $K[t]$ for all $i\le h-1$.
\proof
\begin{enumerate}
\item[$\Leftarrow$)] See Theorem\til\ref{leftarrow}.
\item[$\Rightarrow$)] Suppose $i\le h-2$. By Lemma\til\ref{flatness}, it is enough to show that $E^i_m(M_m)$ is flat over $K[t]/(t^m)$ for all $m\in\Z_+$. We show this by induction on $m\ge 1$:\\
If $m=1$, then $E^i_1(M_1)$ is flat since $K[t]/t$ is a field.\\
Now assume  $m\ge 2$ and  assume that $E^i_{m-1}(M_{m-1})$ is flat over $K[t]/(t^{m-1})$. Since the ideal $(t^{m-1})/(t^m)$ is nilpotent, using the local criterion for flatness (see \cite{matsumura} Theorem $22.3$) we have that $E^i_m(M_m)$ is flat over $K[t]/(t^m)$ if and only if the following two conditions holds:
\begin{enumerate}
\item[$i$)] $E^i_m(M_m)/t^{m-1}E^i_m(M_m)$ is flat over $K[t]/(t^{m-1})$, and
\item[$ii$)] the natural multiplication map
\begin{eqnarray*}
\theta: (t^{m-1})/(t^m)\otimes_{K[t]/(t^m)}E^i_m(M_m)& \la&E^i_m(M_m)\\
\cl{t^{m-1}}\otimes\phi&\map& t^{m-1}\phi
\end{eqnarray*}
is injective, that is, $0:_{E^i_m(M_m)}t^{m-1}=0$.
\end{enumerate}
For each $j\le m-2$, we denote by $\iota^j$ the composition of natural inclusions
$$\iota^j=\iota_j\circ\dots\circ\iota_{m-2}: t^{m-1}M_m\la t^jM_m.$$ 
Since $\mu_j: t^j M_m\la t^{m-1} M_m$ is the multiplication by $t^{m-1-j}$ and $\iota^j$ is the natural inclusion, $\iota^j\circ\mu_j: t^jM_m\la t^jM_m$ is the multiplication by $t^{m-1-j}$.
Furthermore, since $E^i_m(\iota_k)$ is surjective for all $k\in\N$ by Remark\til\ref{rm3} and since $E^i_m(-)$ is a functor, we have that $E^i_m(\iota^j)$ is surjective for all $j\le m-2$. Therefore, since $E^i_m(-)$ is a $A_m$-linear functor, we have
$$\Im(E^i_m(\mu_j))=\Im(E^i_m(\mu_j)\circ E^i_m(\iota^j))=t^{m-1-j}E^i_m(t^jM_m)$$
for all $j\le m-2$. Using again Remark\til\ref{rm3},
$$\Ker(E^i_m(\iota_j))=\Im(E^i_m(\mu_j))=t^{m-1-j}E^i_m(t^jM_m),$$
and so 
$$E^i_m(t^{j+1}M_m)\iso\frac{E^i_m(t^jM_m)}{t^{m-1-j}E^i_m(t^jM_m)}$$
for all $j\le m-2$.\\
Since $m\ge 2$, we can plug in $j=0$ and we get
$$\frac{E^i_m(M_m)}{t^{m-1}E^i_m(M_m)}\iso E^i_m(tM_m).$$
Using Remark\til\ref{rm2} and Remark\til\ref{rm1}, $$E^i_m(tM_m)\iso E^i_m(M_{m-1})\iso E^i_{m-1}(M_{m-1})$$
is flat over $K[t]/(t^{m-1})$ by the inductive hypothesis, and it follows that 
$E^i_m(M_m)/{t^{m-1}E^i_m(M_m)}$ is flat over $K[t]/(t^{m-1})$.
So the condition $i$) is proved. \\
Before proving the condition $ii$), we show first that $\Ker(E^i_m(\iota^j))=tE^i_m(t^jM_m)$ by induction on $j\le m-2$.\\
If $j=m-2$, then $\iota^{m-2}=\iota_{m-2}$ and we have shown above that 
$$\Ker(E^i_m(\iota_{m-2}))=tE^i_m(t^{m-2}M_m).$$
Hence $\Ker(E^i_m(\iota^{m-2}))=tE^i_m(t^{m-2}M_m).$\\
If $j<m-2$, since $E^i_m(\iota^j)=E^i_m(\iota_{m-2})\circ\dots\circ E^i_m(\iota_j)$, we have $$E^i_m(\iota_{m-3})\circ\dots\circ E^i_m(\iota_j)(tx)\in tE^i_m(t^{m-2}M_m)= \Ker(E^i_m(\iota_{m-2}))$$ for each $x\in E^i_m(t^jM_m)$. It follows that $tE^i_m(t^jM_m)\se \Ker(E^i_m(\iota^j))$. On the other hand, if $u\in \Ker(E^i_m(\iota^j))$, then $$E^i_m(\iota_j)(u)\in  \Ker(E^i_m(\iota^{j+1}))=tE^i_m(t^{j+1}M_m)$$
by the inductive hypothesis, and so there exists $v\in E^i_m(t^{j+1}M_m)$ such that $E^i_m(\iota_j)(u)=tv$. Since $E^i_m(\iota_j)$ is surjective by Remark\til\ref{rm3}, there exists $w\in E^i_m(t^{j}M_m)$ such that $E^i_m(\iota_j)(w)=v$. Hence 
$$u-tw\in\Ker(E^i_m(\iota_j))=t^{m-1-j}E^i_m(t^jM_m).$$
It follows that $u\in tE^i_m(t^jM_m)+t^{m-1-j}E^i_m(t^jM_m)=tE^i_m(t^jM_m)$.
Therefore, $tE^i_m(t^jM_m)=\Ker(E^i_m(\iota^j))$ for all $j\le m-2$. \\In particular, we have $tE^i_m(tM_m)=\Ker(E^i_m(\iota^0))$.\\
Now we prove the condition $ii$).
Since 
$$E^i_m(\mu_0)\circ E^i_m(\iota^0)=E^i_m(\mu_0\circ\iota^0): E^i_m(M_m)\la E^i_m(M_m)$$ 
is the multiplication by $t^{m-1}$ and since $E^i_m(\mu_0)$ is injective by Remark\til\ref{rm3}, we have
$$0:_{E^i_m(M_m)}t^{m-1}=\Ker(E^i_m(\mu_0)\circ E^i_m(\iota^0))=\Ker(E^i_m(\iota^0))=tE^i_m(M_m).$$
\endproof
\end{enumerate}
\end{theorem}

\section{Applications}

In this section, we study some applications of $N$-fiber-full modules.

\begin{notations}
Let $R=K[X_1,\dots, X_n]$ be the polynomial ring and fix $w=(w_1,\dots, w_n)\in\N^n$ a weight vector. 
Notice that for each $f\in R$ there exists a unique (finite) subset of the set of monomials of $R$, denoted by $\supp(f)$, such that $$f=\sum_{\mu\in\supp(f)}a_\mu\mu\text{ \quad with\quad }a_\mu\in K\setminus\{0\}.$$
If $\mu=X^u=X_1^{u_1}\cdot\dots\cdot X_n^{u_n}$ with $u=(u_1, \dots, u_n)$, then we set $w(\mu)=w_1u_1+\dots+w_nu_n$. If $f=\sum_{\mu\in\supp(f)}a_\mu\mu\in R$, $f\not=0$, we set $$w(f)=\max\{w(\mu): \mu\in\supp(f)\},$$ 
$$\initi_w(f)=\sum_{\substack{{\mu\in\supp(f)}\\ w(\mu)=w(f)}}a_\mu\mu,$$
and we call 
$$\hom_w(f)=\sum_{\mu\in\supp(f)}a_\mu\mu t^{w(f)-w(\mu)}\in R[t]$$
the \tbs{$w$-homogenization} of $f$.\\
Given an ideal $I\se R$, $\init_w(I)$ denotes the ideal of $R$ generated by $\initi_w(f)$ with  $f\in I$, and $\hom_w(I)$ denotes the ideal of $R[t]$ generated by $\hom_w(f)$ with $f\in I$.
%
\end{notations}
Given a monomial order $<$ on $R=K[X_1,\dots, X_n]$ and given an ideal $I\se R$, there exists a weight vector  $w=(w_1,\dots, w_n)\in(\N_{>0})^n$ such that $\init_<(I)=\init_w(I)$ (see \cite{levico}  \textcolor{black}{Proposition $3.4$}).\\

The following corollary is a generalization of  Corollary $3.3$ in \cite{levico}. 

\begin{corollary}\label{coro1}
Let $I\se R=K[X_1,\dots, X_n]$ be an ideal, $w=(w_1,\dots, w_n)\in\N^n$ a weight vector and suppose that $S=P/\hom_w(I)$ is $N$-fiber-full  up to $h$ as a $P$-module, where $P=R[t]$, $N$ is \fg\ and flat over $K[t]$. Then $\Ext^i_P(S,N)$ is a flat $K[t]$-module for $i\le h-1$ by the previous theorem. So, if furthermore $I$ is homogeneous, we have
$$\dim_K(\Ext^i_R(R/I,N/tN)_j)=\dim_K(\Ext_R^i(R/\init_w(I), N/tN)_j)$$
for all $i\le h-2$ and for all $j\in\Z$.\\
In particular, if $N=P$, then
$$\dim_K(\H_\mm^i(R/I)_j)=\dim_K(\H_\mm^i(R/\init_w(I))_j)$$
for all $i\ge n-h+2$ and for all $j\in\Z$.
\proof
We observe that $\Ext^i_P(S,N)$ is a flat $K[t]$-module for $i\le h-1$ follows directly from Theorem\til\ref{main}. Let us give a graded structure to $R=K[X_1,\dots, X_n]$ by putting $\deg X_i=g_i$ for each $i\in\{1,\dots, n\}$, where $g=(g_1, \dots, g_n)$ is a vector of positive integers. Suppose that $I$ is a $g$-homogeneous ideal,  and note that $\hom_w(I)$ is homogeneous with respect to the bi-graded structure on $P=R[t]$ given by $\deg(X_i)=(g_i, w_i)$ and $\deg(t)=(0,1)$. So $S=P/\hom_w(I)$ and $\Ext_P^i(S, N)$ are \fg\ bi-graded $P$-modules, and it follows that
$$\Ext_P^i(S, N)_{(j, *)}=\bigoplus_{k\in\Z}\Ext_P^i(S, N)_{(j, k)}$$
is a \fg\ graded (with respect to the standard grading) $K[t]$-module for all $j\in \Z$.
By Remark $3.7$ and Remark $3.8$ in \cite{levico},
we have that for each $i,j\in\Z$:
$$\Ext_P^i(S, N)_{(j, *)}\iso K[t]^{a_{i,j}}\oplus(\bigoplus_{k\in\N_{>0}}(K[t]/(t^k))^{b_{i,j,k}})$$
for some natural numbers $a_{i,j}$ and $b_{i,j,k}$. 
Set $b_{i,j}=\sum_{k\in\N_{>0}}b_{i,j,k}$.
Repeating a discussion similar to {Theorem $3.1$} in \cite{levico}
we obtain
$$\dim_K(\Ext^i_R(R/I,N/tN)_j)=a_{i,j}$$
and
$$\dim_K(\Ext_R^i(R/\init_w(I), N/tN)_j)=a_{i,j}+b_{i,j}+b_{i+1,j}$$
for every $i,j\in\Z$.
Since $\Ext^i_P(S,N)$ is a flat $K[t]$-module for $i\le h-1$, $b_{i,j}=0$ for all $i\le h-1$ and so $b_{i+1,j}=0$ for all $i\le h-2$. Hence for all $i\le h-2$ and for all $j\in\Z$
$$\dim_K(\Ext^i_R(R/I,N/tN)_j)=a_{i,j}=\dim_K(\Ext_R^i(R/\init_w(I), N/tN)_j).$$
If $N=P$,  then
$\dim_K(\H_\mm^i(R/I)_j)=\dim_K(\H_\mm^i(R/\init_w(I))_j)$ for all $i\ge n-h+2$ and for all $j\in\Z$
by the local duality theorem for graded modules (see \cite{CMrings} Theorem $3.6.19$).
\endproof
\end{corollary}

In what follows, we suppose furthermore that $ A=\oplus_{i\in\N}A_i$ is positively graded with $A_0=K$ and  $t\in A_1$, and $M$ and $N$ are graded $A$-modules.

\begin{notation}
Let  $\mm$ the homogeneous maximal ideal of $A$, $d=\dim A$ and  let $$I=\bigcap_{i=1}^s\qq_i$$ be a homogeneous ideal of $A$, where $\qq_i$ are the primary components of  $I$. For each integer $h$, we set 
$$I^{\le h}= \bigcap_{\substack{{i=1}\\ \dim(A/\qq_i)\ge d-h}}^s\qq_i.$$
Notice that,
\begin{itemize}
\item $I^{\le d-1}=\sat I$.
\item $\H^i_\mm(A/I)\iso\H^i_\mm(A/I^{\le h})$ for all $i\ge d-h$.
\end{itemize}
\end{notation}

\begin{proposition}\label{prop34}
If $A$ is a $d$-dimensional Cohen-Macaulay ring, $N=\cano_A$ is the canonical module of $A$, $I$ is a homogeneous ideal of $A$ and $A/(I, t)^{\le h}$ is cohomologically full (see Definition $1.1$ in \cite{DDSM}) with $h$ an integer, then $A/I$ is \Nff\ up to $h$.
\proof
First, supposing $m\in\N_{>0}$ and using Theorem $2.3$ in \cite{bounds} we observe that
\begin{eqnarray*}
(I, t^m)^{\le h}&=&\{g\in A\vert \dim\big(A/((I, t^m):g)\big)<d-h\}\\
&=&\{g\in A\vert \dim A-\dim\big(A/((I, t^m):g)\big)>h\}\\
&=&\{g\in A\vert\  \height((I, t^m):g)>h\},
\end{eqnarray*}
in particular,
$$(I, t)^{\le h}=\{g\in A\vert \ \height((I, t):g)>h\}.$$
If $g\in (I, t^m)^{\le h}$, then $\height((I, t^m):g)>h$. Since $(I, t^m)\se (I, t)$, we have $$(I, t^m):g\se (I, t):g,$$ 
it follows that 
$$\height((I, t):g)\ge\height((I, t^m):g)>h.$$
Hence $(I, t^m)^{\le h}\se (I, t)^{\le h}$.\\
We recall that  $f\in\sqrt{(I, t)^{\le h}}$ if and only if there exists $N\in\N$ such that $f^N\in(I, t)^{\le h}$, if and only if  there exists $N\in\N$ such that $\height((I, t):f^N)>h$. Let  $g_1,\dots, g_l$ be the minimal generators of $(I, t):f^N$. We have $g_if^N\in (I, t)$ for all $i\in\{1,\dots, l\}$, hence $$g_i^mf^{Nm}=(g_if^N)^m\in(I, t)^m\se(I, t^m)$$ for all $i\in\{1,\dots, l\}$. It follows that 
$$(g_1^m,\dots, g_l^m)\se (I, t^m):f^{Nm}.$$
Since $$\height((g_1^m,\dots, g_l^m))=\height((I, t):f^N)>h,$$ we have $$\height((I, t^m):f^{Nm})\ge \height((g_1^m,\dots, g_l^m))>h,$$ therefore $f\in\sqrt{(I, t^m)^{\le h}}$. Thus $\sqrt{(I, t)^{\le h}}=\sqrt{(I, t^m)^{\le h}}$.\\
Now for each $m\in\N_{>0}$, we set $\mm_m=\mm/(I,t^m)$.  
Since $A/(I, t)^{\le h}$ is cohomologically full, the natural map $$\H^j_{\mm_m}(A/(I,t^m)^{\le h})\la \H^j_{\mm_1}(A/(I,t)^{\le h})$$
is surjective for all $j$. In general, for all $j$ we have $\H^j_{\mm_1}(A/(I,t))\iso \H^j_{\mm}(A/(I,t))$ and $\H^j_{\mm_m}(A/(I,t^m))\iso \H^j_{\mm}(A/(I,t^m))$, hence the natural map
$$\H^j_{\mm}(A/(I,t^m))\la \H^j_{\mm}(A/(I,t))$$
is surjective for all $j\ge d-h$. Since $A$ is *complete, by the local duality theorem for graded modules (see \cite{CMrings} Theorem $3.6.19$), 
$\Ext^i_A(A/(I,t), N)\la\Ext^i_A(A/(I,t^m), N)$ is injective for all $i\le h$.
\endproof
\end{proposition}

Let $R$  be the polynomial ring $K[X_1,\dots, X_n]$ over a field $K$ and let  $I\se R$ be a homogeneous ideal. In the paper of Conca and Varbaro \cite{CV} they obtained the following result:
\begin{quote}
if $\init(I)$ is a square-free monomial ideal for some term order, then
$\dim_K\H_\mm^i(R/I)_j=\dim_K\H_\mm^i(R/\init(I))_j$ for all $i,j$.
\end{quote}
They actually showed that
$$\dim_K\H_\mm^i(R/I)_j=\dim_K\H_\mm^i(R/\init(I))_j$$ for all $i,j$
if $R/\init(I)$ is cohomologically full. If $\init(I)$ is a square-free monomial ideal then  $R/\init(I)$ is cohomologically full (see \cite{Lyu} Theorem $1$). \\
Using the notion of \Nff\ 
we can prove the following theorem:

\begin{theorem}\label{thm35}
Let $I\se R=K[X_1,\dots, X_n]$ be a homogeneous ideal, $w=(w_1,\dots, w_n)\in\N^n$ a weight vector and  $P=R[t]$. If $R/\init_w(I)^{\le h}$ is cohomologically full, then 
$$\dim_K(\H_\mm^i(R/I)_j)=\dim_K(\H_\mm^i(R/\init_w(I))_j)$$ for all $i>n-h$ and for all $j\in\Z$.\\
In particular, fixing a monomial order $<$ on $R$,
\begin{itemize}
\item if $R/\init_<(I)^{\le h}$ is cohomologically full, then $\dim_K(\H_\mm^i(R/I)_j)=\dim_K(\H_\mm^i(R/\init_<(I))_j)$ for all $i>n-h$ and for all $j\in\Z$;
\item if $\init_<(I)^{\le h}$ is square-free, then $\dim_K(\H_\mm^i(R/I)_j)=\dim_K(\H_\mm^i(R/\init_<(I))_j)$ for all $i>n-h$ and for all $j\in\Z$.
\end{itemize}
\proof
We claim $$R/\init_w(I)^{\le h}\iso P/(\hom_w(I),t)^{\le h+1}.$$
Indeed,  if $\init_w(I)=\bigcap_{i=1}^s\qq_i$ is a primary decomposition of  $\init_w(I)$, by definition
$$\init_w(I)^{\le h}= \bigcap_{\substack{{i=1}\\ \dim(R/\qq_i)\ge n-h}}^s\qq_i.$$
Since $\qq_i\se R$, $x, t$ is a $P$-regular sequence for each $x\in R$ such that $x\not=0$ and $x$ is not invertible,  we have  that $(\qq_i, t)$ are  the primary components of the ideal $(\init_w(I), t)$ and 
\begin{eqnarray*}
(\init_w(I)^{\le h}, t)&=&\big(\bigcap_{\substack{{i=1}\\ \dim(R/\qq_i)\ge n-h}}^s\qq_i\big)+(t)\\
&=&\bigcap_{\substack{{i=1}\\ \altezza(\qq_i)\le h}}^s(\qq_i,t)\\
&=&\bigcap_{\substack{{i=1}\\  \altezza(\qq_i, t)\le h+1}}^s(\qq_i,t)\\
&=&(\init_w(I), t)^{\le h+1}.
\end{eqnarray*}
By Proposition $3.5$ in \cite{levico}, we have $$(\hom_w(I),t)=(\init_w(I),t).$$ 
Hence
$$R/\init_w(I)^{\le h}\iso P/(\init_w(I)^{\le h},t)\iso P/ (\init_w(I), t)^{\le h+1}\iso P/(\hom_w(I),t)^{\le h+1}.$$
Since $R/\init_w(I)^{\le h}$ is cohomologically full, by our claim and by Proposition\til\ref{prop34} we have that $P/\hom_w(I)$ is $P$-fiber-full up to $h+1$. Therefore 
$$\dim_K(\H_\mm^i(R/I)_j)=\dim_K(\H_\mm^i(R/\init_w(I))_j)$$
for all $i> n-h$ and for all $j\in\Z$ by Corollary\til\ref{coro1}.\\
In particular, given a monomial order $<$ on $R$, there exists a weight vector $w=(w_1,\dots, w_n)\in(\N_{>0})^n$  such that $\init_w(I)=\init_<(I)$, and this yields the statement.
\endproof
\end{theorem}

\begin{remark}
In the paper of  Dao, De Stefani and Ma  \cite{DDSM} they proved the following result (see \cite{DDSM} Lemma $3.7$):
\begin{quote}
Let $R$  be the polynomial ring $K[X_1,\dots, X_n]$ over a field $K$ and let  $J\se R$ be a homogeneous ideal. If $R/J$ is a cohomologically full ring, then $R/J$ satisfies Serre's condition ($S_1$), that is, $\min(J)=\ass(R/J)$.
\end{quote}
Hence, it can happen that $R/J$ is not  cohomologically full only because it has some embedded primes. But using the notion $J^{\le h}$, sometimes we can remove the embedded primes of  $R/J$ and make  $R/J^{\le h}$ be a cohomologically full ring.
\end{remark}

\begin{example}
Let $I$ be a monomial ideal of the polynomial ring $R=K[x_1,\dots,x_n ]$  and let $\mu_1, \dots, \mu_s$ be  the minimal monomial generators of $I$ such that the following condition holds:
if there exists $\mu\in\{\mu_1,\dots, \mu_s\}$ and there exists an integer $t\ge 2$ such that $x_k^t\vert\mu$ with $k\in\{1,\dots, n\}$, then 
there exists $g\in\sat I$ such that $g\vert\sqrt{\mu}$. Then  we have $\sat I=I^{\le n-1}$ is square-free. Hence  $R/\sat I$ is cohomologically full.\\
A simple example is the following: If $R=K[x,y,z]$ and $J=(x^2y, xy^2,xyz)$, then we have $\ass(R/J)=\{(x, y, z), (x), (y)\}$ and $\min(J)=\{(x), (y)\}$, hence $R/J$ is not a  cohomologically full ring. We observe $\ass(R/\sat J)=\{(x), (y)\}=\min(\sat J)$. Since $\sat J$ is square-free, $R/\sat J$ is cohomologically full.
\end{example}

\begin{remark}
We observe that in the situation of Corollary\til\ref{coro1}, if $N=K[t]$ then 
$$\dim_K(\Ext^i_R(R/I,N/tN)_j)=\dim_K(\Ext^i_R(R/I,K)_j)=\beta_{i,j}(R/I)$$
is the $(i,j)$-th Betti number of $R/I$ and 
$$\dim_K(\Ext_R^i(R/\init_w(I), N/tN)_j) = \beta_{i,j}(R/\init_w(I))$$
is the $(i,j)$-th Betti number of $R/\init_w(I)$, hence
$$\beta_{i,j}(R/I)= \beta_{i,j}(R/\init_w(I))$$ for all $i\le h-2$ and for all $j\in\Z$. 
However in general, $\beta_{i,j}(R/I)\not= \beta_{i,j}(R/\init_w(I))$ even if $\init_w(I)$ is a square-free monomial ideal.
Hence $R$ is fiber-full does not imply $R$ is $K[t]$-fiber-full. \\
In practice, if $A$ is Cohen-Macaulay and $A/(I, t)$ is cohomologically full, then the natural map
$\Ext^i_A(A/(I,t), \cano_A)\la\Ext^i_A(A/J, \cano_A)$ is injective for all $i$,
where $J\se(I, t)$ and $\sqrt{J}=\sqrt{(I,t)}$. Hence, if $(I,t)$ is a square-free monomial ideal, then the natural map
$\Ext^i_A(A/(I,t), N)\la\Ext^i_A(A/J, N)$ is injective if $N=\cano_A$, however it is not true if $N=K[t]$. 
\end{remark}

%
%

\begin{example}
Consider the following graph $G$:
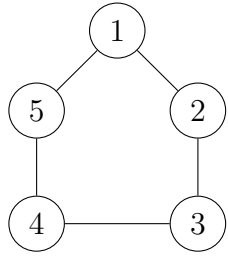
\begin{center}
\begin{tikzpicture}[node distance =1.5 cm, main/.style = {draw, circle}] 
\node[main] (1) {$1$}; 
\node[main] (2) [below right of=1] {$2$};
\node[main] (3) [below of=2] {$3$}; 
\node[main] (5) [below left of=1] {$5$}; 
\node[main] (4) [below of=5] {$4$};
\path (1) edge [bend right =0] node[above] {$ $} (5);
\path (1) edge [bend left =0] node[above] {$ $} (2);
\path (5) edge [bend right =0] node[above] {$ $} (4);
\path (2) edge [bend left =0] node[above] {$ $} (3);
\path (4) edge [bend right =0] node[above] {$ $} (3);
\end{tikzpicture} 
\end{center}
We obtain the binomial edge ideal of $G$ (see \cite{bino}):
$$J_G=(x_1y_2-x_2y_1, x_2y_3-x_3y_2, x_3y_4-x_4y_3,  x_4y_5-x_5y_4, x_1y_5-x_5y_1),$$
and using Macaulay2 \cite{M2} we compute the initial ideal of $J_G$:
\begin{eqnarray*}
\init_<(J_G) &=&(x_2y_1, x_5y_1, x_3y_2, x_1x_5y_2, x_4y_3, x_1x_2x_5y_3, \\
&&x_5y_4, x_4y_1y_5, x_1x_4y_2y_5, x_3y_1y_4y_5),
\end{eqnarray*}
where $<$ is  the lexicographic order on $K[x_1, x_2,\dots, x_5, y_1, y_2, \dots, y_5]$ induced by $$x_1>x_2>\dots>x_5>y_1>y_2>\dots>y_5.$$
Therefore, $\init_<(J_G)$  is a square-free monomial ideal, 
$\beta_0(J_G)=5$ and $\beta_0(\init_<(J_G))=10\not=\beta_0(J_G)$. 
\end{example}

Hence, in the situation of Corollary\til\ref{coro1}, it would be very interesting to understand when  $S=P/\hom_w(I)$ is $K[t]$-fiber-full, a condition that would guarantee that the graded Betti numbers are preserved going from $I$ to $\init_w(I)$.

\end{document}